\documentclass[11pt]{amsart}
\usepackage{amssymb,latexsym}

\newcommand{\iso}{{\;\;\stackrel{_\sim}{\longrightarrow}\;\;}}
\newcommand{\vi}{${\sf {(i)}}\;$}
\newcommand{\vii}{${\sf {(ii)}}\;$}
\newcommand{\viii}{${\sf {(iii)}}\;$}

\newcommand{\too}{\,\,\longrightarrow\,\,}
\newcommand{\onto}{\,\,\twoheadrightarrow\,\,}
\newcommand{\reg}{^{^{\sf{reg}}}}

\newcommand{\BA}{{\mathbb A}}
\newcommand{\BC}{{\mathbb C}}
\newcommand{\BN}{{\mathbb N}}

\newcommand{\BP}{{\mathbb P}}
\newcommand{\BZ}{{\mathbb Z}}

\newcommand{\CC}{{\mathcal C}}
\newcommand{\CCM}{{\mathcal {CM}}}

\newcommand{\CG}{{\mathcal G}}

\newcommand{\CO}{{\mathcal O}}

\newcommand{\irrep}{{\sf{Irrep}}}
\newcommand{\Gr}{{\mathsf {Gr}}}

\newcommand{\Mat}{{\mathsf {Mat}}}
\newcommand{\Sch}{{\mathsf {Sch}}}

\newcommand{\GL}{{\mathbf {GL}}}
\newcommand{\PGL}{{\mathbf {PGL}}}

\newcommand{\fC}{{\overline{\mathcal C}}{}}
\newcommand{\fh}{{\mathfrak h}}
\newcommand{\fS}{{\mathfrak S}}
\newcommand{\fP}{{\mathfrak P}}

\newcommand{\fm}{{\mathfrak m}}

\newcommand{\oA}{{\overset{_\circ}{\mathbb A}}{}}

\newcommand{\oSch}{{\overline{\mathsf {Sch}}}{}}

\begin{document}

\title[]{Calogero-Moser space and Kostka polynomials}
\author{Michael Finkelberg}
\address{Independent Moscow University, Bolshoj Vlasjevskij pereulok, dom 11,
Moscow 121002 Russia}
\email{fnklberg@mccme.ru}
\author{Victor Ginzburg}
\address{Department of Mathematics, University of Chicago, Chicago, IL 60637,
USA}
\email{ginzburg@math.uchicago.edu}

\maketitle
\begin{abstract}We consider 
the canonical map from the Calogero-Moser
space to symmetric powers of the affine line, sending
conjugacy classes of pairs of $n\times n$-matrices
to their eigenvalues. 
We show that the character of a 
natural $\BC^*$-action on the scheme-theoretic
zero fiber of this  map is given by
Kostka polynomials.
\end{abstract}

\section{Introduction}

\subsection{}
The aim of this paper is to prove a refined version of Conjecture 17.14
of ~\cite{eg}. To explain our result, recall the
so-called {\it Calogero-Moser space} ${\mathcal C}_n$, a
$2n$-dimensional
complex
algebraic manifold introduced by Kazhdan-Kostant-Sternberg ~\cite{kks},
and studied further by G. Wilson ~\cite{w}. 
It is defined as 
 $\CC_n:=\CCM_n/\!/\PGL_n$, the  quotient by
the natural (free) conjugation-action of the group
$\PGL_n$ on the set
$$\CCM_n:= 
\lbrace{(X,Y)\in \Mat_n\times\Mat_n\quad\big|\quad
[X,Y]+{\sf{Id}}=\;\text{rank } 1 \text{ matrix}\rbrace}.
\eqno(1.1)$$

Let
 $\BA^{(n)}$ denote the set of unordered
$n$-tuples of complex numbers.
The assignment $(X,Y)\mapsto \bigl(Spec(X)\,,\,Spec(Y)\bigr)$,
sending a pair of $n\times n$-matrices  to the
corresponding pair of
$n$-tuples of their eigenvalues
gives a map $p: \CC_n \to \BA^{(n)}\times \BA^{(n)}$.
The zero fiber $p^{-1}(0,0)$ of this map
is formed by the conjugacy classes of {\it nilpotent}
pairs $(X,Y)\in \CCM_n$. This fiber is known to be
a finite set labelled naturally by partitions of
$n$. Given such a partition $\lambda$,
let  $p^{-1}(0,0)_\lambda$ be the corresponding
point in the  zero fiber.

The Conjecture 17.14 of ~\cite{eg} states that,
for any  partition $\lambda$,
the corresponding
point in the (scheme-theoretic) zero fiber of $p$ comes with
{\it multiplicity} equal to $(\dim V_\lambda)^2$,
where $V_\lambda$ is an irreducible representation
of the symmetric group $\fS_n$ attached 
to the  partition $\lambda$ in the standard way, see ~\cite{m}.

\subsection{}
In the present paper we propose and prove the following
$q$-analogue of the above conjecture.
Observe that the complex torus $\BC^*$ acts naturally
on $\CCM_n$ by $z: (X,Y) \mapsto (z^{-1}\cdot X\,,\,z\cdot Y)
\,,\,\forall z\in \BC^*.$
This $\BC^*$-action descends to the 
Calogero-Moser space ${\mathcal C}_n$
and preserves the  zero fiber $p^{-1}(0,0)$.
Now given $\lambda$, a partition of $n$,
let  $p^{-1}(0,0)_\lambda$ be the corresponding
irreducible component of the  zero fiber
viewed as a {\it non-reduced} scheme
(set theoretically concentrated at one point).
The $\BC^*$-action keeps theese points (set theoretically) fixed
hence, for each $\lambda$, induces a $\BC^*$-action on the coordinate ring
of the scheme
$p^{-1}(0,0)_\lambda$, a finite dimensional vector space.
The character of this finite dimensional $\BC^*$-module
may be viewed as a Laurent polynomial $ch_\lambda\in
\BZ[q,q^{-1}]$.
Now, recall that for each partition $\lambda$ 
one defines the {\em Kostka polynomial} $K_\lambda(q)\in \BZ[q]$
which is a certain 
$q$-analogue\footnote{In the main body of the paper we use a minor
modification of the standard $K_\lambda(q)$.} 
of $\dim V_\lambda$, the dimension  of the corresponding irreducible
$\fS_n$-representation, see e.g. (\cite{m}~III.6).  

Our result reads

\smallskip {\bf Theorem.} {\em For any partition $\lambda$ (of $n$),
we have: $ch_\lambda=K_\lambda(q)\cdot K_\lambda(q^{-1})$.}

\subsection{}
This result has a natural generalization to
other finite complex reflection groups $W$
in a vector space $\fh$. In more details, in ~\cite{eg}
the authors associate to a pair $(\fh, W)$
a {\em Calogero-Moser space}
$\CC_W$ together with a finite map $p:\ \CC_W\to\fh/W\times\fh/W$.
In the special case $\fh=\BC^n$ and
$W=\fS_n,$ the space $\CC_W$ reduces to the
variety $\CC_n$,
and the map $p:\CC_W\to\fh/W\times\fh/W$
reduces to the map $p:  \CC_n \to \BA^{(n)}\times \BA^{(n)}$
considered above.

 More generally, in this paper we will consider the case where
$\fh=\BC^n$ and $W=\Gamma\sim\fS_n$ 
is a wreath product of $\fS_n$ and $\Gamma=\BZ/N\BZ$,
a cyclic group of some fixed order $N$
(thus, $W=\fS_n\ltimes(\BZ/N\BZ)^n$),
acting naturally in $\fh$.
It has been proved in ~\cite{eg} that the
corresponding Calogero-Moser space
$\CC_{\Gamma,n}:=\CC_W$ is
a smooth affine algebraic variety isomorphic to a
certain Nakajima's Quiver variety for a cyclic quiver. 

The Conjecture 17.14 of ~\cite{eg} says that the reduced fiber of
$p$ over $(0,0)\in\fh/W\times\fh/W$ can be identified with the set
$\irrep(W)$ of isomorphism classes of irreducible representations of $W$,
and the multiplicity of the point in this fiber
corresponding to $\rho\in
\irrep(W)$ 
equals $(\dim \rho)^2$. 

It is well known that the irreducible $\Gamma\sim\fS_n$-modules are
naturally parametrized by the set $\fP_\Gamma(n)$ of 
$\Gamma^\vee$-partitions of $n$, see e.g. (\cite{m},
Part ~I, Appendix B). Here $\Gamma^\vee$ is the set of irreducible 
characters of $\Gamma$, and a $\Gamma^\vee$-partition $\Lambda$ is a
collection $(\lambda_\chi,\ \chi\in\Gamma^\vee)$ of ordinary partitions
such that $\sum_\chi|\lambda_\chi|=n$. It is known that 
the points of reduced fiber of $\CC_{\Gamma,n}$ over $(0,0)$ are
also naturally numbered by $\fP_\Gamma(n)$ (in case of trivial $\Gamma$
it was proved in ~\cite{w}, and in the general case in ~\cite{k}). 
By abuse of notation we will denote the
point in the fiber corresponding to $\Lambda\in\fP_\Gamma(n)$ by $\Lambda$
as well. 

\subsection{}
\label{main}
The cyclic Calogero-Moser space $\CC_{\Gamma,n}$ has
a natural $\BC^*$-action, such that its
 fixed point set $\CC_{_{\Gamma,n}}^{^{\BC^*}}$
coincides with the reduced zero fiber. We consider
the character of induced $\BC^*$-action in the Artin coordinate ring
$\CO_\Lambda$ of the component $p^{-1}(0,0)_\Lambda$
of the fiber concentrated
at the point $\Lambda\in\CC_{_{\Gamma,n}}^{^{\BC^*}}$.

\subsection{}
\label{second}
For an arbitrary cyclic group $\Gamma=\BZ/N\BZ$ and
$\Lambda\in\fP_\Gamma(n)$, 
we introduce a polynomial
$K_\Lambda(q)\in\BZ[q]$
which is a $q$-analogue of $\dim V_\Lambda$,
 the corresponding
irreducible $\Gamma\sim\fS_n$-module, see ~\ref{analog}, 
and prove the following

\smallskip {\bf Theorem.} {\em The character of $\BC^*$-module
$\CO_\Lambda$ equals $K_\Lambda(q)\cdot K_\Lambda(q^{-1})$.}

\subsection{} Our proof is a straightforward application of the remarkable
work ~\cite{w}. G.~Wilson has studied the reduced 
fibers of the second projection
$p_2:\ \CC_n\to\BA^{(n)}$ and identified them  
as certain products of Schubert cells in Grassmannians. His results
reduce our problem to some classical computations in Grassmannians.

One ingredient in the proof of Theorem \ref{main}
is a relative {\em Drinfeld
compactification} $\fC_n$ of the Calogero-Moser space $\CC_n$
(such that
the projection $p_2$ extends to the proper projection 
$p_2:\ \fC_n\to\BA^{(n)}$, see ~\ref{drinf}) and its cyclic version,
see ~\ref{cyclic Dr}. Though it enters our proof only at some technical
point, we believe that $\fC_n$ is a very interesting object in itself.

The space  $\fC_n$ was, in fact,
 implicitly introduced in ~\cite{w} where Wilson studied
the embedding of $\CC_n$ into the {\em adelic Grassmannian} $\Gr_{ad}$
(the cyclic version of this embedding is studied in ~\cite{bgk}).
Wilson constructed a  set-theoretic partition 
$\Gr_{ad}=\bigsqcup_{k\in\BN}\CC_n$. However,
it turns out that the union
$\bigsqcup_{0\leq k\leq n}\CC_k$ {\em can not} be equipped with the
 structure
of an
algebraic variety. The algebraic variety $\fC_n$ has,
on the other hand, a natural partition
$\fC_n=\bigsqcup_{0\leq k\leq n}\CC_k\times\BA^{(n-k)},$ 
see ~\ref{twist}, into smooth locally-closed strata
(similar in spirit to the stratification used in  \cite{ku})
and may be viewed as an algebraic `resolution' of 
$\bigsqcup_{0\leq k\leq n}\CC_k$,
a nonalgebraic substack of $\Gr_{ad}$. The name "{\it
Drinfeld's compactification}"
is suggested by a close analogy with Drinfeld's quasimap spaces,
cf. ~\cite{ku}. 

In ~\ref{hope} we propose an alternative
 conjectural definition of $\fC_n$ as a step
towards its generalization for other
 Nakajima quiver varieties.

\subsection{Acknowledgments} We are deeply obliged to 
A.~Kuznetsov for very useful discussions and comments.
M.F. is grateful to the University of Chicago for the wonderful
working conditions, and to V.~Vologodsky for patient explanations
of the trivia of intersection theory. This research was conducted by M.F.
for the Clay Mathematics Institute.

\section{Wilson's embedding into a relative Grassmannian}

\subsection{The Calogero-Moser space} 
\label{CM}
 Fix  a positive integer $n$ and consider
the space $\CCM_n$ defined in (1.1).
 Then $\CCM_n$ is smooth, and
the action of $\PGL_n$ by the simultaneous 
conjugation is free (see ~\cite{w}).
The quotient space $\CC_n:=\CCM_n/\PGL_n$ is a $2n$-dimensional 
smooth affine algebraic variety,
the {\em Calogero-Moser space}. For $n=0$ we define $\CC_0$ to be a point.

Recall that $\BA^{(n)}:=\BA^n/\fS_n$.
The assignment $Y\mapsto Spec(Y)$,
sending a matrix $Y\in \Mat_n$ to the
$n$-tuple of its eigenvalues viewed
as a finite subscheme 
 of $\BA^1$ given by zeros of the characteristic
polynomial of $Y$, yields
an isomorphism of algebraic varieties:
$\Mat_n//\PGL_n \iso\BA^{(n)}$ 
(where $\Mat_n//\PGL_n$ denotes 
the categorical quotient). The
second projection $\CCM_n\to\Mat_n\,,\, (X,Y) \mapsto Y,$ 
descends to the projection
$\pi_n:\ \CC_n\to\BA^{(n)}$. Wilson has
determined all the reduced fibers of $\pi_n$. 
Namely, he constructed an embedding of
any fiber into a certain product of (finite dimensional) Grassmann varieties,
and identified the image with a union of products of certain Schubert cells.
Let us formulate his results more precisely.
Till the end of this section 
{\em fiber}
means {\em reduced fiber}, and we write $\pi^{-1}(-)$
instead of $\pi^{-1}(-)_{_{\sf{reduced}}}$.

\subsection{Theorem} 
\label{factor}
(Wilson, ~\cite{w}, 7.1)
{\em Suppose a divisor $D=D_1+D_2\in\BA^{(n)}$ 
is a sum of divisors $D_1\in\BA^{(m)},\ D_2\in\BA^{(k)}$
with
disjoint supports. Then there is a
canonical isomorphism $\pi_n^{-1}(D)\simeq
\pi_m^{-1}(D_1)\times\pi_k^{-1}(D_2)$.}

\medskip

We will refer to this result
as  the {\em factorization property}
of the projection $\pi_n$ (or rather of
the collection of maps $\pi_n$ over
$n\in\BN$).

\subsection{}
\label{fiber}
In view of the above theorem, in order to describe an arbitrary fiber of
$\pi_n$, it suffices to describe the fiber over the principal diagonal,
$\pi_n^{-1}(ny),\ y\in\BA^1$. To this end, consider the
polynomial algebra $\BC[z]$ and,
for any $y\in\BC$ write $\fm_y=(z-y)\cdot\BC[z] $ for the corresponding
maximal ideal. Let $\Gr(n,y)$ be
Grassmannian  of $n$-dimensional subspaces in the
vector  space
$\BC[z]/\fm_y^{\,2n}$. The  vector space $\BC[z]/\fm_y^{\,2n}$
comes equipped with a distinguished
complete flag
$$0\subset\fm_y^{\,2n-1}/\fm_y^{\,2n}\subset
\fm_y^{\,2n-2}/\fm_y^{\,2n}\subset\ldots\subset
{\fm_y}/\fm_y^{\,2n}\subset\BC[z]/\fm_y^{\,2n}
$$
(quotients of ideals). This flag defines the {\em Schubert stratification}
of $\Gr(n,y)$. Let $\Sch_n(y)\subset\Gr(n,y)$ denote the
locally closed subvariety formed by all the Schubert cells of dimension $n$.

\smallskip 

{\bf Theorem.} (\cite{w}, 6.4) {\em There is a
canonical isomorphism $\pi_n^{-1}(ny)\simeq\Sch_n(y)$.}

\subsection{}
\label{adel}
Wilson also describes the way the above fibers glue together. In order
to formulate his result, we recall that $\BA^{(n)}$ may be viewed as the
space of all codimension $n$ ideals $I\subset\BC[z]$, and introduce the
following definition\smallskip

{\bf Definition.} The {\em relative Grassmannian} $\CG_n$ is the space
of pairs $(I,W)$ where $I\subset\BC[z]$ is a codimension $n$ ideal, and
$W\subset\BC[z]/I^2$ is an $n$-dimensional linear subspace.

\medskip

Clearly, $\CG_n$ is a quasiprojective variety equipped with a projection
$\pi_n:
\CG_n\onto\BA^{(n)}\,,\,(I,W)\mapsto I.$ For any 
$I\in\BA^{(n)}$ we have:
 $\pi_n^{-1}(I)\simeq\Gr(n,2n)$.

Wilson considers an open subset $\CC_n\reg\subset\CC_n$ formed by the
(conjugacy classes of) pairs $(X,Y)$ such that $Y$ is diagonalizable and has
pairwise
distinct eigenvalues. Each element in $\CC_n\reg$ 
has a unique representative
of the form $Y=diag(y_1,\ldots,y_n),\ X=\|x_{ij}\|,$ with
$x_{ij}=(y_i-y_j)^{-1},$ for $i\ne j,$ and
$x_{ii}=\alpha_i$. 
To the $n$-tuple $(y_1,\ldots,y_n)$ we
associate the  $n$-tuple
of lines $W_i=\big\langle 1-\alpha_i(z-y_i)\big\rangle
\subset \BC[z]/\fm_{y_n}^{\,2}\,,\,i=1,\ldots,n,$
in the corresponding 2-planes.
Wilson defines an embedding $\beta:\ \CC_n\reg\to\CG_n$
by the formula $\beta: (X,Y)\mapsto (I,W)$ where 
$I=(z-y_1)\cdot\ldots\cdot(z-y_n)$,
and $W$ is set to be a direct sum of the lines $W_i$, that is:
$$I=\fm_{y_1}\cdot\ldots\cdot\fm_{y_n}\enspace,\enspace
W=\oplus_i\,W_i\,\,
\subset\BC[z]/I^2\equiv\BC[z]/\fm_{y_1}^{\,2}\oplus\ldots\oplus
\BC[z]/\fm_{y_n}^{\,2}.$$

\smallskip {\bf Theorem.} (Wilson, ~\cite{w}, 5.1) 
\vi {\it The map $\beta$ extends to 
an embedding $\beta:\ \CC_n\hookrightarrow\CG_n$ commuting with the projections
$\pi_n$;}

\vii {\it Given  $D=\sum_{k=1}^l n_ky_k\in\BA^{(n)}$
and $C\in\pi_n^{-1}(D)\subset\CC_n$
write $\ C=(W_1,\ldots,W_l),\ 
W_k\in\Sch_{n_k}(y_k)$. Then, 
 under the natural identification
$\BC[z]/\prod_{k=1}^l\fm_{y_k}^{\,2n_k}\,\equiv\,\bigoplus_{k=1}^l
\BC[z]/\fm_{y_k}^{\,2n_k}$, we have}
$$\beta:\; C\;\mapsto \;\bigl(\prod\nolimits_{k=1}^l\,\fm_{y_k}^{\,n_k}\;\;,\;\;
\bigoplus\nolimits_{k=1}^l\,\,W_k\bigr)\,.\hspace{3cm}\Box$$

\subsection{Drinfeld relative compactification}
\label{drinf}
We define $\fC_n\subset\CG_n$ as the closure of $\beta(\CC_n)$ or,
equivalently, of $\beta(\CC_n\reg)$. 
Specifically, consider the open stratum of the diagonal stratification
$\oA^{(n)}\subset\BA^{(n)}$ formed by all the $n$-tuples of 
pairwise distinct points.
Consider the locally closed subvariety $\fC_n\reg\subset\pi_n^{-1}(\oA^{(n)})
\subset\CG_n$ formed by all the pairs
$$
\Big\lbrace\!\displaystyle{\begin{array}{ll}\displaystyle
\bigl(I,W)& \,\big|_{_{}}\;
I=\fm_{y_{1_{}}}\cdot\ldots\cdot\fm_{y_n}\quad,\quad W
\subset\BC[z]/I^2\equiv_{_{}}\\
&\BC^{}[z]/\fm_{y_1}^{\,2^{}}\oplus\ldots\oplus
\BC[z]/\fm_{y_n}^{\,2^{}}\bigr),\enspace
\text{such that}\enspace
W\cap (\BC[z]/\fm_{y_i}^{\,2})\neq 0\,,\,\forall i.
\end{array}}\!\Big\rbrace
$$
Thus, $W$ is a direct
sum of lines $W_i\subset\BC[z]/\fm_{y_i}^{\,2}$.
\smallskip

{\bf Definition.} The {\em Drinfeld compactification} 
$\fC_n\subset\CG_n$ is defined as the closure of $\fC_n\reg$ in $\CG_n$.
The restriction of $\pi_n:\ \CG_n\to\BA^{(n)}$ to $\fC_n$ is also 
denoted by $\pi_n$. 
\smallskip

Clearly, $\pi_n:\ \fC_n\to\BA^{(n)}$ is a projective morphism.

\subsection{Twist by a divisor}
\label{twist}
The rest of this 
section will not be used elsewhere in
 the paper but it helps to 
understand better the structure of $\fC_n$.

For $0\leq k\leq n$ we will define a map 
${\mathsf{twist}}_{k}^n:\ 
\fC_k\times\BA^{(n-k)}\to\fC_n$ (twist by a divisor).
To this end,
given an ideal $I\subset\BC[z]$ of codimension $n-k$, and
$(J,W)\in\fC_k$,
take the preimage of $W$ under the natural projection
$\BC[z]/IJ^2\to\BC[z]/J^2$, and let
$W'\subset I/I^2J^2
\subset\BC[z]/I^2J^2$ correspond to this preimage under the natural
identification $I/I^2J^2\simeq\BC[z]/IJ^2$. 
We set ${\mathsf{twist}}_k\bigl((J,W),I\bigr):=(IJ,W')$.

    From now on we will 
identify $\CC_n$ with its image $\beta(\CC_n)\subset
\fC_n\subset\CG_n$. Given $y\in\BA$,
write $\oSch_m(y)\subset\Gr(m,\BC[z]/\fm_y^{\,2m})$
for the union of Schubert cells of dimension $\leq m$.
Wilson's theorem ~\ref{adel} immediately implies
the following
\smallskip 

{\bf Theorem.} \vi {\em Let $D=\sum_{k=1}^l n_ky_k\in\BA^{(n)}$. Then
$\pi_n^{-1}(D)\subset\fC_n$ equals $\prod_{k=1}^l\oSch_{n_k}(y_k)$.
 Specifically, under the natural identification 
$\BC[z]/\prod_{k=1}^l\fm_{y_k}^{2n_k}\equiv\bigoplus_{k=1}^l
\BC[z]/\fm_{y_k}^{2n_k}$, a point $(W_1,\ldots,W_l),\
W_k\in\oSch_{n_k}(y_k)$, corresponds to 
$\bigoplus_{k=1}^l\,W_k$.

\vii $\fC_n\smallsetminus \CC_n=
{\mathsf{twist}}_{n-1}^n(\fC_{n-1}\times\BA^1)$,
where the RHS is a closed subvariety.

\viii $\fC_n$ is a disjoint union of the locally closed subvarieties:}
\[\fC_n=\bigsqcup\nolimits_{k=1}^n\;
{\mathsf{twist}}_{k}^n(\CC_k\times\BA^{(n-k)})\,.\]

Part (i) implies, in particular, that the map $\pi_n:\ \fC_n\to\BA^{(n)}$
enjoys the factorization property.

\subsection{Remark}
\label{hope}
One would like to find a construction of $\fC_n$ in the ordinary
Calogero-Moser setup of ~\ref{CM}, avoiding the use of adelic Grassmannian.
Here is a conjectural definition. Recall that 
$\CCM_n\subset\Mat_n\times\Mat_n$ is a smooth closed subvariety.
Now $\Mat_n$ can be viewed
as an open subset of $\Gr(n,2n)$ via identifying a matrix $X$
with the graph
$W_X\subset\BC^n\oplus\BC^n$ 
of the corresponding linear map $\BC^n\to\BC^n$. Let $\CCM'_n$ be the 
closure of $\CCM_n$ in $\Gr(n,2n)\times\Mat_n$. The group $\GL_n$ acts on
$\CCM'_n$ naturally: $g(W,Y)=(gW,gYg^{-1})$. Let $\fC'_n$ be the GIT
quotient of $\CCM'_n$ with respect to $\GL_n$.
\smallskip

{\bf Question.} {\em Is there an isomorphism $\fC'_n\simeq\fC_n$ 
which is the identity on
the common open subset $\fC'_n\supset\CC_n\subset\fC_n$?}
 
\section{$\BC^*$-action on Schubert cells}

\subsection{} Write $\Gr(n)$ instead of 
 $\Gr(n,0)$ for the Grassmannian of $n$-dimensional
subspaces of $\BC[z]/(z^{2n})$. We have the 
standard complete flag in $\BC[z]/(z^{2n})$ (see ~\ref{fiber}):
$$0\subset\fm_y^{\,2n-1}/\fm_y^{\,2n}\subset
\fm_y^{\,2n-2}/\fm_y^{\,2n}\subset\ldots\subset
{\fm_y}/\fm_y^{\,2n}\subset\BC[z]/\fm_y^{\,2n}
$$
Recall that 
$\Sch_n\subset\Gr(n)$ is a disjoint union of the $n$-dimensional cells,
which are known to be exactly the cells $\Sch_\lambda$ numbered
by the set $\fP(n)$ of partitions of $n$. In more detail, 
given a partition $\lambda=(l_1,\ldots,
l_n),\ 0\leq l_1\leq\ldots\leq l_n,\ l_1+\ldots+l_n=n$, we have
$$\Sch_\lambda=\{
W\in\Gr(n)\quad\mid\quad
\dim\bigl(W\cap(\fm_0^{2n-l_i-i}/\fm_0^{2n})\bigr)=i\,,\; \forall i
=1,\ldots,n\}. 
$$
The multiplicative group $\BC^*$ acts on $\BC[z]$ by $(c,z^i)\mapsto
c^{-i}z^i$. This action induces a natural action on $\Sch_\lambda
\subset\Gr(n)$ contracting this Schubert cell to the unique fixed point
$W_\lambda:=\langle z^{2n-l_1-1},z^{2n-l_2-2},\ldots,z^{2n-l_n-n}\rangle$.
The tangent space $T_{W_\lambda}\Sch_\lambda$ at the
point $W_\lambda$ is naturally isomorphic to

\begin{align*}
&Hom(\BC z^{2n-l_1-1},\langle z^{2n-1},\ldots,z^{2n-l_1}\rangle)\times\\
&Hom(\BC z^{2n-l_2-2},
\langle z^{2n-1},\ldots,\widehat{z^{2n-l_1-1}},\ldots,z^{2n-l_2-1}\rangle)
\times\ldots\times\\
&Hom(\BC z^{2n-l_n-n},
\langle z^{2n-1},\ldots,\widehat{z^{2n-l_1-1}},\ldots,\widehat{z^{2n-l_i-i}},
\ldots,z^{2n-l_n-n+1}\rangle)
\end{align*}
where $\;\widehat{\vphantom{x}}\;$ means omission of an element.
   From this we read off easily the character of $\BC^*$ on 
$T_{W_\lambda}\Sch_\lambda$. Specifically, 
write $h_\lambda(u)$ for the hook length of a
box $u$ in the Young diagram attached naturally to a partition
$\lambda$.
Below, we use the notation $q^i$
for the character $\BC^*\to\BC^*\,,\,c\mapsto c^i,$
and write $ch^{\,}V$ for the character of a finite dimensional
$\BC^*$-module $V$.

\subsection{Lemma}
\label{hooks} We have:
$ch^{\,}(T_{W_\lambda}\Sch_\lambda)=
\sum_{u\in\lambda}\,\,q^{-h_\lambda(u)}.\hfill\Box$

\section{Nilpotent extensions of Schubert cells}

\subsection{} 
\label{projections}
Recall the map $\pi_n:\ \CC_n\to\BA^{(n)}\,,\,
(X,Y)\mapsto Spec(Y).$ Denote this map by $p_2$, and
similarly, consider the other
projection
$p_1:\ \CC_n\to\BA^{(n)}\,,\,
(X,Y)\mapsto Spec(X).$ 
 Note that there is an involution $\omega$ on $\CC_n$
such that $\omega:
(X,Y)\mapsto (Y^t,X^t),$ and we have: $p_1=p_2\circ\omega$. Let $p=(p_1,p_2)$
stand for the simultaneous
projection  $(p_1,p_2): \CC_n\too\BA^{(n)}\times\BA^{(n)}$.
To distinguish between the two copies of $\BA^{(n)}$ we will use the
notation $p:\ \CC_n\to\BA^{(n)}_1\times\BA^{(n)}_2$. 
According to ~\cite{eg}, the map $p$ is a finite morphism. 

The scheme theoretic fiber $p_2^{-1}(0)$ is a disjoint union of schemes
$p_2^{-1}(0)_\lambda$ such that the underlying reduced scheme is 
$\Sch_\lambda$, to be denoted $\Sch_\lambda^2$ from now on. Similarly, 
the scheme theoretic fiber $p_1^{-1}(0)$ is a disjoint union of schemes
$p_1^{-1}(0)_\lambda$ such that the underlying reduced scheme is 
denoted by $\Sch_\lambda^2$. 

Our goal is to compute the scheme theoretic fiber $p^{-1}(0,0)$. It is well
known that the corresponding reduced scheme is a disjoint union of points:
the $\BC^*$-fixed points of $\Sch_n^1$ (or equivalently, $\Sch_n^2$). 
Abusing the language we will denote the $\BC^*$-fixed point of 
$\Sch_\lambda^2$ by $\lambda$; thus $\Sch_\lambda^1\cap\Sch_\lambda^2=
\lambda$. We will denote the connected component of $p^{-1}(0,0)$
concentrated at $\lambda$ by $p^{-1}(0,0)_\lambda$. 

Note that $p^{-1}(0,0)_\lambda$ is the fiber over $0\in\BA^{(n)}_1$ with
respect to the projection $p_1:\ p_2^{-1}(0)_\lambda\to\BA^{(n)}_1$.
Our first step will be to compute the fiber over $0\in\BA^{(n)}_1$ with
respect to the projection $p_1:\ \Sch_\lambda^2\to\BA^{(n)}_1$.
 
\subsection{}
\label{kost}
Recall that the Kostka polynomial associated to a Young diagram
$\lambda$ is a polynomial in the variable `$q$' given by the formula:
$\displaystyle{q^{m(\lambda)}(1-q)
\ldots(1-q^n)\prod\nolimits_{u\in\lambda}(1-q^{h_\lambda(u)})^{-1}},
$
where  $m(\lambda)$ is
a certain positive integer, see (\cite{m}, page 243, Example 2). 
This is a  $q$-analogue
of the dimension $\dim 
V_\lambda$ of the irreducible representation $V_\lambda$
of the symmetric group $\fS_n$. We will consider a version of Kostka
polynomial with the lowest term equal to 1, that is, we put
$$K_\lambda(q):=(1-q)\ldots(1-q^n)
\prod\nolimits_{u\in\lambda}(1-q^{h_\lambda(u)})^{-1}\,.
$$

{\bf Proposition.} {\em We have:}
$$ch^{\,}\CO(\Sch_\lambda^2\cap p_1^{-1}(0)_\lambda)=
K_\lambda(q)\quad\text{and}\quad
ch^{\,}\CO(\Sch_\lambda^1\cap p_2^{-1}(0)_\lambda)
=K_\lambda(q^{-1})\,.$$

{\em Proof.} The two formulas are analogous, so we only prove the
first one. We  compute the geometric fiber  of the sheaf
$(p_{1})_*\CO(\Sch_\lambda^2)$ at the point $0\in\BA^{(n)}_1$. This is a 
locally free coherent sheaf, that is a (trivial) vector bundle, so to 
compute the character of its geometric fiber  
at 0 it suffices to know the character $ch^{\,}\CO(\Sch_\lambda^2)$
of its space of global sections,
and the character of $\CO(\BA^{(n)}_2)$. 
Now we  pass to the formal completions at 0 and $\lambda$.
Thus, we are reduced to finding
 the characters of tangent spaces $T_0\BA^{(n)}_1$
and $T_\lambda\Sch_\lambda^2$. The former character equals 
$1+q^{-1}+\ldots+q^{-n}$, while the latter character was computed in
the Lemma ~\ref{hooks}. We conclude that 
$ch_{\,}\hat\CO_{\BA^{(n)}_1,0}=(1-q)^{-1}\ldots(1-q^n)^{-1},$ and
$ch_{\,}\hat\CO_{\Sch_\lambda^2,\lambda}=
\prod_{u\in\lambda}(1-q^{h_\lambda(u)})^{-1}$.
Thus, we get 
\[ch_{\,}\CO(\Sch_\lambda^2\cap p_1^{-1}(0)_\lambda)=
ch_{\,}\hat\CO_{\Sch_\lambda^2,\lambda}/ch_{\,}\hat\CO_{\BA^{(n)}_1,0}=
K_\lambda(q)\,.\qquad \Box\]

\subsection{} 
\label{kost2}
We are going to compute $ch_{\,}\CO(p^{-1}(0,0)_\lambda)$ along
similar lines. To this end it suffices to compute the character of
the completion $ch_{\,}
\hat\CO_{p_2^{-1}(0)_\lambda,\lambda}$. We will prove that
$ch_{\,}\hat\CO_{p_2^{-1}(0)_\lambda,\lambda}=
K_\lambda(q^{-1})\prod_{u\in\lambda}(1-q^{h_\lambda(u)})^{-1}$. Hence, 
arguing exactly as in the proof of ~\ref{kost} we
will be able to conclude that
$ch_{\,}\CO(p^{-1}(0,0)_\lambda)=K_\lambda(q)K_\lambda(q^{-1})$, as required
in the Theorem ~\ref{main}. Thus, to prove the Theorem it suffices
to prove the following
\medskip

{\bf Proposition.} $ch_{\,}\hat\CO_{p_2^{-1}(0)_\lambda,\lambda}=
K_\lambda(q^{-1})\prod_{u\in\lambda}(1-q^{h_\lambda(u)})^{-1}$.

\subsection{} 
\label{generic}
We start the proof the Proposition with the following 
\smallskip

{\bf Lemma.} {\em $\Sch_\lambda^1$ and $\Sch_\lambda^2$ are transversal at
$\lambda$.}
\smallskip

{\em Proof.} The varieties
$\Sch_\lambda^1$ and $\Sch_\lambda^2$ 
are smooth of complementary dimensions. Moreover, the character
of $T_\lambda\Sch_\lambda^2$ is a polynomial in $q^{-1}$ without constant
term, while the character
of $T_\lambda\Sch_\lambda^1$ is a polynomial in $q$ without constant
term. Hence, these two tangent spaces must have zero intersection,
and we are done. \hfill$\Box$\smallskip

Thus the formal completion of $\CC_n$ at $\lambda$ is isomorphic to
a product of formal
completions of $\Sch_\lambda^1$ and $\Sch_\lambda^2$ at $\lambda$.
We will denote by $pr_1$ and $pr_2$ the projections arising this way. 
The fiber over $\lambda$ of the
restriction of $pr_2$ to the formal completion of $p_2^{-1}(0)_\lambda$
equals:
$pr_2^{-1}(\lambda)=\Sch_\lambda^1\cap p_2^{-1}(0)_\lambda$.
 We already know formulas for
$ch_{\,}\CO(\Sch_\lambda^1\cap p_2^{-1}(0)_\lambda)$ and 
$ch_{\,}\hat\CO_{\Sch_\lambda^2,\lambda}$, so to complete the proof it suffices
to show that $pr_{2*}\hat\CO_{p_2^{-1}(0)_\lambda,\lambda}$ is a (trivial)
vector bundle on the completion of $\Sch_\lambda^2$ at $\lambda$. To this
end it suffices to show that the dimension of the generic 
 fiber  of
$(pr_{2})_*\hat\CO_{p_2^{-1}(0)_\lambda,\lambda}$ equals\linebreak
$\dim\CO(\Sch_\lambda^1\cap p_2^{-1}(0)_\lambda)=d_\lambda\,
(=\dim V_\lambda).$ But the 
dimension of the generic  fiber  equals $m_\lambda$,
the multiplicity of the scheme
$p_2^{-1}(0)$ at the generic point of its reduced 
subscheme $\Sch_\lambda^2$.

To compute this multiplicity $m_\lambda$ we may as well work in the
Drinfeld compactification $\fC_n$ embedded into the relative Grassmannian
$\CG_n$ over $\BA^{(n)}_2$. A general fiber $p_2^{-1}(\underline{y})
\subset\Gr(n,\underline{y})$
is reduced at the generic point, so $m_\lambda$ is  the coefficient of
the cycle class $[p_2^{-1}(\underline{y})]$ 
with respect to the Schubert basis 
$\{[\oSch_\lambda],\ \lambda\in\fP(n)\}$ of the degree $2n$
homology group of $\Gr_n$.

Now recall that a general $n$-tuple
$\underline{y}=(y_1,\ldots,y_n)\in\BA^{(n)}_2$ 
of pairwise distinct points gives rise to a
direct sum  decomposition
$\BC[z]/\fm_{y_1}^{\,2}\cdot\ldots\cdot\fm_{y_n}^{\,2}$
$=\bigoplus_i\,\BC[z]/\fm_{y_i}^{\,2}$, and
$p_2^{-1}(\underline{y})
\subset\Gr(n,\underline{y})$ is  the product of corresponding projective
lines: $p_2^{-1}(\underline{y})=
\BP^1\times\ldots\times\BP^1\subset\Gr_n$. It is the
classical result of Schubert calculus that 
for the corresponding homology classes one has
an expansion:
$[\BP^1\times\ldots\times\BP^1]=\sum_\lambda m_\lambda\cdot
[\oSch_\lambda],$
where the coefficients $m_\lambda$ 
can be read off from the formula:
 ${\mathbf p}_1^n= \sum_\lambda \,m_\lambda\cdot s_\lambda,$
an expansion
of the $n$-th power of the first symmetric function
 ${\mathbf p}_1^n$
with respect to the basis of Schur functions $s_\lambda$. 
The  coefficients in the the latter  expansion 
are well-known to be
 equal to $d_\lambda=K_\lambda(1)$, see e.g. (\cite{m}, page 114).

This completes the proof of Proposition ~\ref{kost2} and the proof of
Theorem ~\ref{main}. \hfill$\Box$\break

\section{Cyclic Calogero-Moser space}

\subsection{} 
\label{cyclic CM}
Consider the action of $\Gamma=\BZ/N\BZ\subset\BC^*$ on
the Calogero-Moser space $\CC_{nN}$. The fixed-point subvariety 
$\CC_{nN}^\Gamma$ consists of various connected components. 
There is a single component 
characterized by the property that the representation of $\Gamma$
in the fiber of tautological bundle 
at any point in this component is a multiple of the
regular representation, see ~\cite{k}. We will call this connected component
$\CC_{\Gamma,n}$. According to {\em loc. cit.}, $\CC_{\Gamma,n}$ is a special
case of Nakajima's Quiver variety (corresponding to $N$-cyclic quiver with
$n$-dimensional spaces at all "finite" vertices,  1-dimensional space at
an "extended" 
vertex, and a nonzero value of the diagonal moment map). 

We have the natural projection $p=(p_1,p_2):\ \CC_{\Gamma,n}\to
(\BA^{(nN)}_1\times\BA^{(nN)}_2)^\Gamma$. Note that 
$$(\BA^{(nN)}_1\times\BA^{(nN)}_2)^\Gamma=
(\BA^{(nN)}_1)^\Gamma\times(\BA^{(nN)}_2)^\Gamma
\quad\text{and}\quad
(\BA^{(nN)}_{1,2})^\Gamma=(\BA^1_{1,2}/\Gamma)^{(n)}\,.
$$
We let $\BA_{\Gamma,1,2}^{(n)}$ denote
the set on the right of this formula,  and view $p$ as a projection 
$p=(p_1,p_2):\ \CC_{\Gamma,n}\to\BA^{(n)}_{\Gamma,1}\times
\BA^{(n)}_{\Gamma,2}$. The natural $\BC^*$-action  on $\CC_{nN}$
when restricted to $\CC_{\Gamma,n}$ factors through 
$\BC^*\stackrel{c\mapsto c^N}{\longrightarrow}\BC^*$, and we will
consider the resulting $\BC^*$-action on $\CC_{\Gamma,n}$ (which is 
generically free).

\subsection{}
\label{cyclic Gr}
Wilson's embedding $\beta:\ \CC_{nN}\hookrightarrow\CG_{nN}$ is 
$\Gamma$-equivariant, and its image lands into a connected component
$\CG_{\Gamma,n}\subset\CG_{nN}^\Gamma$ 
characterized by the property that the representation of $\Gamma$
in the fiber of tautological bundle 
at any point of this component is a multiple of the
regular representation (to see the inclusion:
 $\CC_{\Gamma,n}\subset\CG_{\Gamma,n}$
it suffices to check it at any $\BC^*$-fixed point, e.g. $\lambda=(nN)$). 
We will denote by $\beta:\ \CC_{\Gamma,n}
\hookrightarrow\CG_{\Gamma,n}$ this cyclic version of Wilson's embedding,
and we will use it to describe the reduced fibers of $p_2$.

First of all, the action of $\Gamma$ on $\BC[z]/(z^{2nN})$
splits it into a direct sum $\BC[z]/(z^{2nN})=
\bigoplus_{\chi\in\Gamma^\vee}\BC[z]/(z^{2nN})_\chi$
of $N$ $2n$-dimensional eigenspaces according
to the characters of $\Gamma$. Note that we can
canonically identify $\Gamma^\vee$ with $\BZ/N\BZ$, and then 
$\BC[z]/(z^{2nN})_\chi$ is spanned by $\{z^k,\ k\equiv\chi\pmod{N}\}$.
The fiber of $\CG_{\Gamma,n}$ over $nN\cdot0\in(\BA_2^{(nN)})^\Gamma$
is formed by all the $nN$-dimensional subspaces $W\subset\BC[z]/(z^{2nN})$
such that $W=\bigoplus_{\chi\in\Gamma^\vee}W_\chi,\ W_\chi\subset
\BC[z]/(z^{2nN})_\chi$, and $\dim W_\chi=n$. Thus, this fiber equals
$\prod_{\chi\in\Gamma^\vee}\Gr(n,\BC[z]/(z^{2nN})_\chi)$. 
Each space $\BC[z]/(z^{2nN})_\chi$ carries a natural complete flag
(given
by the
intersections with  powers of the maximal ideal). Thus, each
variety
$\Gr(n,\BC[z]/(z^{2nN})_\chi)$ has a natural stratification into Schubert
cells numbered by partitions. 

Set $\fP_\Gamma(n):=\{\lambda_\chi\;\mid\;
\chi\in\Gamma^\vee\,,\,\sum_\chi|\lambda_\chi|=n\}$, and given
 $\Lambda\in\fP_\Gamma(n)$ put
$\Sch_\Lambda^2:=
\prod_\chi\Sch_{\lambda_\chi}\subset
\prod_{\chi\in\Gamma^\vee}\Gr(n,\BC[z]/(z^{2nN})_\chi).$
Now Wilson's theorem ~\ref{fiber} 
together with ~\cite{k} yield the following.
\smallskip

{\bf Proposition.} {\em The reduced fiber of $\CC_{\Gamma,n}$ over
$0\in\BA_{\Gamma,2}^{(n)}$ is canonically isomorphic to 
$\coprod_{\Lambda\in\fP_\Gamma(n)}\;
\Sch_\Lambda^2.\hfill\Box$}
\smallskip

{\bf Corollary.} \vi
{\em Each component $\Sch_\Lambda^2$ contains a unique
$\BC^*$-fixed point $\Lambda\in\CC_{\Gamma,n}$.} 

\vii {\em The reduced fiber
of $p:\ \CC_{\Gamma,n}\to\BA_{\Gamma,1}^{(n)}\times\BA_{\Gamma,2}^{(n)}$
over $(0,0)$ coincides with $\CC_{_{\Gamma,n}}^{^{\BC^*}}
=\fP_\Gamma(n).\hfill\Box$}
\smallskip

We will denote by $p^{-1}(0,0)_\Lambda$ the connected component of the 
scheme theoretic fiber concentrated at the point $\Lambda$, and we will
denote by $p_{1,2}^{-1}(0)_\Lambda$ the connected component of the
scheme theoretic fiber concentrated at $\Sch_\Lambda^{1,2}$.

\subsection{}
\label{cyclic Dr} 
We define the Drinfeld compactification 
$\fC_{\Gamma,n}\supset\CC_{\Gamma,n}$ as the closure of $\CC_{\Gamma,n}$
inside $\CG_{\Gamma,n}$. 

We will need a description of a general fiber of $p_2:\ \fC_{\Gamma,n}\to
\BA_{\Gamma,2}^{(n)}$. If we choose a primitive $N$-th root of unity $\zeta$
then a general point $\underline{y}\in\BA_{\Gamma,2}^{(n)}$ 
can be represented by a
collection 
$$\underline{y}=(y_1,\zeta y_1,\ldots,\zeta^{N-1}y_1,y_2,\ldots,\zeta^{N-1}y_2,
\ldots,y_n,\ldots,\zeta^{N-1}y_n)$$ of distinct points of $\BA_2^1$.
The $2nN$-dimensional vector space 
$V=\BC[z]/\fm_{y_1}^{\,2}\ldots\fm_{\zeta^{N-1}y_n}^{\,2}$ 
is acted upon by $\Gamma$,
and splits up into a direct sum of $N$ $2n$-dimensional eigenspaces $V_\chi$
according to the characters of $\Gamma$. We also have a direct sum
decomposition $V=U_1\oplus\ldots\oplus U_n$ where 
$U_i=\BC[z]/\fm_{y_i}^{\,2}\ldots
\fm_{\zeta^{N-1}y_i}^{\,2}$. Note that for any $i$ and 
$\chi$ the intersection $U_i\cap V_\chi$ is 2-dimensional. We will denote
this intersection by $UV_{i,\chi}$.

The fiber of $\CG_{\Gamma,n}$ over $\underline{y}$ equals 
$\prod_{\chi\in\Gamma^\vee}\Gr(n,V_\chi)$. The fiber of $\fC_n$ over 
$\underline{y}$ is isomorphic to  $\prod_{1\leq i\leq n}\BP^1$; let us
explain how it is embedded into $\prod_{\chi\in\Gamma^\vee}\Gr(n,V_\chi)$.
We have a direct sum decomposition $U_i=\bigoplus_{k\in\BZ/N\BZ}
\BC[z]/\fm_{\zeta^ky_i}^{\,2}$, and the action of $\Gamma$ on $U_i$ permutes
these summands. Hence $\BC[z]/\fm_{y_i}^{\,2}$ projects isomorphically onto any
$UV_{i,\chi}$. Given a line $\ell_i\in\BP^1(\BC[z]/\fm_{y_i}^{\,2})$ we denote
by $\ell_{i,\chi}\subset UV_{i,\chi}$ its image under the above isomorphic
projection. Finally, for a collection $(\ell_i)\in\prod_{1\leq i\leq n}
\BP^1(\BC[z]/\fm_{y_i}^{\,2})$ the corresponding point of 
$\prod_{\chi\in\Gamma^\vee}\Gr(n,V_\chi)$ is the collection of subspaces 
$(\bigoplus_i\ell_{i,\chi}\subset V_\chi)$.

\subsection{}
\label{analog}
Our aim is to compute the character of $\BC^*$-action on the Artin ring
$\CO(p^{-1}(0,0)_\Lambda)$, that is, to prove Theorem ~\ref{second}.
The proof is entirely similar to that of ~\ref{main}. Let us  spell out
the intermediate steps. 
First, we define: 
$$K_\Lambda(q):=
(1-q)\ldots(1-q^n)
\prod\nolimits_{\chi\in\Gamma^\vee}^{u\in\lambda_\chi}
\;(1-q^{h_\lambda(u)})^{-1}\,.
$$
Analoguously to Proposition ~\ref{kost} we obtain
\smallskip

{\bf Proposition.} {\em We have:}
$$
ch^{\,}\CO(\Sch^2_\Lambda\cap p_1^{-1}(0)_\Lambda)=
K_\Lambda(q)
\quad\text{and}\quad
ch^{\,}\CO(\Sch^1_\Lambda\cap p_2^{-1}(0)_\Lambda)
=K_\Lambda(q^{-1}).\hfill\Box
$$

Further, an analogue of Proposition ~\ref{kost2} reads

\subsection{Proposition}
\label{last}
 $\displaystyle{ch_{\,}\hat\CO_{p_2^{-1}(0)_\Lambda,\Lambda}=
K_\Lambda(q^{-1})
\prod\nolimits_{\chi\in\Gamma^\vee}^{u\in\lambda_\chi}
\,(1-q^{h_\lambda(u)})^{-1}\,.}$
\smallskip

To prove this last Proposition we argue, as in ~\ref{generic}, that it suffices
to check if the generic multiplicity $m_\Lambda$ of $p_2^{-1}(0)_\Lambda$
equals $d_\Lambda:=K_\Lambda(1)$. To this end we turn to the cyclic version
of Drinfeld compactification $\fC_{\Gamma,n}$, see ~\ref{cyclic Dr}. 
A general fiber $p_2^{-1}(\underline{y})$ being reduced at the generic point,
$m_\Lambda$ are the coefficients of the cycle class 
$[p_2^{-1}(\underline{y})]$ with respect to the Schubert basis
$\{[\oSch_\Lambda],\ \Lambda\in\fP_\Gamma(n)\}$ 
of the degree $2n$  homology group of 
$\prod_{\chi\in\Gamma^\vee}\Gr_n$. Our description of the general fiber
$p_2^{-1}(\underline{y})$ in ~\ref{cyclic Dr} boils down to the following.

Take the diagonal embedding $\BP^1=\Delta_{\BP^1}\hookrightarrow
\prod_{\chi\in\Gamma^\vee}\BP^1_\chi$. For each $\chi\in\Gamma^\vee$ we have
an  embedding $(\BP^1_\chi)^n\hookrightarrow\Gr_n$ as in ~\ref{generic}.
Now form the composition 
$$(\BP^1)^n=(\Delta_{\BP^1})^n\;\hookrightarrow\;
\prod\nolimits_{\chi\in\Gamma^\vee}\,(\BP^1_\chi)^n\;\hookrightarrow\;
\prod\nolimits_{\chi\in\Gamma^\vee}\,\Gr(n,V_\chi)\,.
$$ 

The homology class of $[\Delta_{\BP^1}]$ in the 2-homology of 
$\prod_{\chi\in\Gamma^\vee}\BP^1_\chi$ equals $\sum_\chi[\BP^1_\chi]$, the sum
of degree 2 the  generators of the
 homology groups of the factors. As in ~\ref{generic}, 
we conclude that $[(\Delta_{\BP^1})^n]=\sum_\Lambda m_\Lambda
\cdot[\oSch_\Lambda]$
where the coefficients $m_\Lambda$ equal the coefficients of 
$(\sum_\chi {\mathbf p}_{1,\chi})^n$ 
with respect to the basis of Schur functions
$S_\Lambda$ (here ${\mathbf p}_{1,\chi}$ 
is the first power sum symmetric function
in the variables $x_{i,\chi},\ 1\leq i<\infty$, and 
$S_\Lambda=\prod_\chi s_{\lambda_\chi}(x_{i,\chi})$, see ~\cite{m}, part I,
Appendix B.) The latter coefficients are in turn equal to:
$n!/\prod_{\chi\in\Gamma^\vee}^{u\in\lambda_\chi}h_\lambda(u)
=d_\Lambda=K_\Lambda(1)$,
see {\em loc. cit.} (9.6) on page 178.

This completes the proof of Proposition ~\ref{last}, hence the proof
of
Theorem ~\ref{second}. $\Box$

\end{document}